\documentclass[a4paper,11pt]{article}
 
\usepackage{color}

\newcommand{\ZD}{\;\mbox{\rm d}}    
\usepackage[latin1,utf8]{inputenc}
\usepackage[english]{babel}
\usepackage{amsthm,amssymb,amsmath,amsopn,amsfonts}
\theoremstyle{plain}
\newtheorem{theorem}{Theorem}[]

\theoremstyle{definition}

\theoremstyle{remark}

\newcommand{\F}{\Phi}

\title{%
{\bf Regularity of the steering control for systems with persistent memory}
\thanks{This papers fits into the research program of the GNAMPA-INDAM and has been written in the framework of the   ``Groupement de Recherche en Contr\^ole des EDP entre la France et l'Italie (CONEDP-CNRS)''.}%
} 
\author{%
{L. Pandolfi} \thanks{Politecnico di Torino, Dipartimento di matematica, Corso Duca degli Abruzzi~24, 10129 Torino Italy,
  {\tt luciano.pandolfi@polito.it}}
\and
{D. Triulzi}\thanks{Graduated from the Dipartimento di matematica ``Federigo Enriquez'', Via saldini 50, 20133 Milano, {\tt daniele.triulzi1@studenti.unimi.it }}
}  

\date{~}
\begin{document}
\maketitle

\begin{abstract}
The following fact is known for large classes of distributed control systems: when the target is regular, there exists a regular steering control. This fact is important to prove convergence estimates of numerical algorithms for the approximate computation of the steering control.

We   extend this property to a class of systems with persistent memory  of Maxwell-Boltzmann type. 
\end{abstract}

\section{Introduction}

We study the following system where $x\in (0,\pi)$ and $t>0$:

\begin{equation}\label{equazione}
\begin{cases}
w''(x,t)=w_{xx}(x,t)+\int_0^tM(t-s)w_{xx}(x,s)\ZD s\,,\quad w(0,t)=f(t)\,, \quad w(\pi,t)=0 \\
w(x,0)=0\,, \quad w'(x,0)=0\,.\\
\end{cases}
\end{equation}
We assume $M(t)\in H^2(0,T)$ and  $f(t)\in L^2(0,T)$ for every $T>0$.   As proved for example in~\cite{Pandolfi}, 
  $w(x,t)\in C([0,T];L^2(0,\pi))\cap C^1([0,T];H^{-1}(0,\pi))$ and 
  for every $(\xi,\eta)  \in L^2(0,\pi)\times  H^{-1}(0,\pi)$ and $T>2\pi$, there exists $f\in L^2(0,T)$ such that  $
w(T)=\xi\,,\quad w'(T)=\eta 
$. We prove:

\begin{theorem}\label{teo:regularity}
Let $(\xi,\eta)\in H^1_0(0,\pi)\times L^2(0,\pi)$ and let $T>2\pi $. There exists a steering control $f\in H^1_0(0,T)$.

\end{theorem}

 The proof is in Section~\ref{Section:momentREGULARITY}. Extensions to the case $x\in\Omega\subseteq \mathbb{R}^d$, $d>1$, is reserved to future investigations.

 We conclude this introduction with few comments. First we note that system~(\ref{equazione}) is often encountered in the study of viscoelasticity and diffusion equations with memory. When $M(t)=0$ of course it reduces to the string equation. In the case of the wave equation (even when $x$ in regions of  $\mathbb{R}^d$, $d>1$) theorem~\ref{teo:regularity} is known. The proof that we give, based on moment methods, shows in particular controllability (in $H^1_0(0,\pi)\times L^2(0,\pi)$) of the cascade connection of system~(\ref{equazione}) with an integrator. We refer to~\cite[Ch.~11]{Tucsnak} and references therein for these facts.

\section{\label{Section:momentREGULARITY}The   proof of Theorem~\ref{teo:regularity}}

The following computations are a bit simplified if 
we integrate the first equation of (\ref{equazione}) on $[0,t]$
and we write it in the equivalent form  (here $N(t)=1+\int_0^tM(s)\ZD s$)

\begin{equation}\label{eq:formaPrimoordine}
w'(x,t)=\int_0^tN(t-s)w_{xx}(x,s)\ZD s\,, \qquad 
w(x,0)=0\,,\quad 
w(0,t)=f(t)\,,\  w(\pi,t)=0 \,.
\end{equation}

We use  the orthonormal basis of $L^2(0,\pi)$ whose elements are $
\F_n=\sqrt{(2/\pi)}\sin nx$, $ n\in \mathbb{N}$, and we expand

\[
w(x,t)=\sum_{n\in\mathbb{N}}\F_n(x)w_n(t)\,, \quad 
w_n(t) =\sqrt{\frac{2}{\pi}}
\int_0^\pi\F_n(x)w(x) \ZD x\,. 
\]
Then $w_n(x,t)$ must satisfy

\[    
w'_n(t) =-n^2\int_0^tN(t-s)w_n(s)\ZD s 
 +n\int_0^tN(t-s)\left (\sqrt{2/\pi}f(s)\right )\ZD s \,.
\]
The function $ \sqrt{2/\pi}f  $ will be renamed $f$. 

Let  $z_n(t)$ solve
\begin{equation}\label{eq:dizN} 
z'_n(t)=-n^2\int_0^tN(t-s)z_n(s)\ZD s\,, \quad 
z_n(0)=1\,.
 \end{equation}
We have  (see~\cite{PandIEOT})

\begin{align} 
\nonumber w_n(t)&=n\int_0^t\left(\int_0^{t-s}N(t-s-\tau)z_n(\tau)d\tau\right)f( s)\ZD s=\\
\label{wn}&= \frac{1}{n}\int_0^t\left ( \frac{\ZD}{\ZD s}z_n (t-s)\right ) f(s) \ZD s\,,
\\
\label{w'n}
w_n'(t)&=n\int_0^t\left(-\frac{\ZD}{ \ZD s} \int_0^{t-s}N(t-s-\tau)z_n(\tau)d\tau \right)
f(s)\ZD s\,.
\end{align}

We require that  a target $(\xi,\eta)\in H^1_0(0,\pi)\times L^2(0,\pi)$  is reached at time $T$, i.e. we require   $(w(T),w'(T))=(\xi,\eta)$. 

The Fourier expansion of the targets is 

  \[
  \xi=\sum_{n=1}^{+\infty}\frac{\xi_n}{n}\F_n\, ,\quad  \mbox{and}\quad \eta=\sum_{n=1}^{+\infty} \eta_n \F_n\,,
\qquad    (\{\xi_n\},\{\eta_n\})\in l^2(\mathbb{N})\times l^2(\mathbb{N} )
  \,. 
  \] 
 So, controllability to $(\xi,\eta)$ at time $T$ is equivalent to the existence   of a control $f\in L^2(0,T)$ such that $w_n(T)=\xi_n/n$, $w'_n(T)=\eta_n$ for every $n$. The expression we found for $w_n(t)$ and $w_n'(t)$ suggest that we investigate whether is it possible to solve this problem with
   \begin{equation}\label{eq:ContH1}
   f(t)=\int_0^tg(s)\ZD s\,,\quad g\in L^2(0,T)\,. 
    \end{equation}
 If this is possible then we have the existence of an $H^1$-steering control, and we get a steering control in $H^1_0(0,T)$ if we can find $g$ which satisfies the additional condition
 \begin{equation} \label{eq:medianulla}
 \int_0^T g(s) \ZD s=0\,.
 \end{equation}
We replace the expression~(\ref{eq:ContH1}) in $w_n(T)$ and $w_n'(T)$ and we integrate by parts. We see that  $f$ is an $H^1$ steering control to $(\xi,\eta)$ if the following
\emph{moment problem} is solvable: 
 
\begin{align}
 \xi_n & 
  =\int_0^T g(r)\ZD r-\int_0^T z_n(T-s)g(s)\ZD s\,,
\label{xi}\\
   \eta_n&= \int_0^T\left [n\int_0^{T-s} N(T-s-r) z_n(r)\ZD r\right ] g(s)\ZD s=
\label{eta}  
 \int_0^Tg(T-s)\left (\frac{-z_n'(s)}{n} \right )\ZD s\,  .
\end{align}

We multiply equation (\ref{eta}) by $i$ and we subtract it from   (\ref{xi}). Furthermore we impose the additional condition~(\ref{eq:medianulla}). We find the moment problem:

\begin{equation}\label{momenti}
\int_0^TZ_n(s) g(T-s)\ZD s= c_0\,,\quad c_0=\left\{\begin{array}{lll}
\xi_n-i\eta_n&{\rm if}& n>0\\
0&{\rm if}& n=0
\end{array}\right.
\end{equation} 
and $Z_n(t)=\left(z_n(s)+\frac{i}{n}z'_n(s)\right)$ if $n>0$, $Z_0(t)=1$\,. 
In order to prove statement~\ref{teo:regularity:ITEM1} of Theorem~\ref{teo:regularity}, we prove solvability of the moment problem~(\ref{momenti}).

We note that $\{c_n\}_{n>0}$ is an arbitray \emph{complex valued} $l^2(\mathbb{N})$ sequence while $g$ is real (when $\xi$ and $\eta$ are real). We reformulate the moment problem (\ref{momenti}) with $n\in \mathbb{Z} $. We proceed as follows: for $n< 0$ we define:

\[
z_n(t)=z_{-n}(t),\quad \F_n(x)=\F_{-n}(x)\,,\quad Z_{-n}(t)=\bar{Z}_n(t)\,.
\]
 
Proceeding as in~\cite[Lemma~5.1]{Pandolfi} we can see that the moment problem~(\ref{momenti}) can be equivalently studied with $n\in\mathbb{Z}$ and $g$ complex valued.

Our goal is the proof that the moment problem~(\ref{momenti}),    $n\in\mathbb{Z}$, is solvable.
Even more, we prove that  
$\{ Z_n(t)\}_{n\in \mathbb{Z}}$
is a Riesz sequence in $L^2(0,T)$, provided that $T>2\pi$. This  imply the   additional information that \emph{the solution $g\in L^2(0,T)$ of minimal norm  depends continuously on the target $(\xi,\eta)\in H^1_0(0,\pi)\times L^2(0,\pi)$. }
  Integrating this function $g$ as in~(\ref{eq:ContH1}) we get the steering control $f$ of minimal norm in $H^1_0(0,T)$ and so \emph{the solution $f\in H^1_0(0,T)$ of minimal norm  depends continuously on the target $(\xi,\eta)\in H^1_0(0,\pi)\times L^2(0,\pi)$. }

\subsection{The proof that $\{Z_n\}_{n\in\mathbb Z}$ is a Riesz sequence in $L^2(0,T)$, $T>2\pi$}

The proof that $\{Z_n\}_{n\in\mathbb Z}$ is a Riesz sequence in $L^2(0,T)$, $T>2\pi$,   is divided in two steps: in the first one we show that the sequence $\{Z_n\}_{n\in\mathbb{Z}\setminus\{0\}}$ is a Riesz sequence in $L^2(0,T)$. Then we will prove that $\{ Z_n\}_{n\in\mathbb{Z} }$ is a Riesz sequence in $L^2(0,T)$ too.  
In the proof we use the  following definitions and results (see~\cite[Chp.~3]{Pandolfi}):  
a sequence $\{x_n\}$ in a Hilbert space $H$ is: 
\begin{itemize}
\item
a \emph{Riesz sequence} when it is the image of an orthonormal sequence under a linear bounded and boundedly invertible transformation;
\item
 $\omega$-\emph{independent} when the following holds: if $\{ \alpha_n \} \in l^2$ and if $\sum_{n=1}^{+\infty}\alpha_n x_n=0$ (convergence in the norm of $H$) then   $\{\alpha_n\}=0$.  
 \end{itemize}
\emph{Let $\{x_n\}$ be a Riesz sequence in the Hilbert space $H$ and let $\{y_n\}$ be quadratically close to $\{x_n\}$, i.e. $\sum_{n=1}^{+\infty}\|x_n-y_n\|^2_{H} <+\infty$ . Then 
  there exists $N$ such that  $\{y_n\}_{|n|>N}$ is a Riesz sequence. 
 If furthermore  $\{y_n\}$ is $\omega$-independent   then it is a Riesz sequence too.}

We introduce the notation
and 
$
\mathbb{Z}'=\mathbb{Z}\setminus\{0\} 
$. 
 
\paragraph{Step 1: $\{Z_n\}_{n\in \mathbb{Z}'}$ is a Riesz sequence in $L^2(0,T)$, $T>2\pi$  }  This part of the proof is contained in~\cite{Pandcina}. For completeness, we report the proof  in the form we need here.

We put $N'(0)=\gamma$.  Using~\cite[Lemmas~5.2 and~5.5]{Pandcina} we get that for every $T>0$ there exists $C$ such that
\begin{equation}\label{eq:DiseqAsint}
\sum_{n\in\mathbb{Z}'}\left\|Z_n(t)-e^{\gamma t}e^{int}\right\|^2_{L^2(0,T)}\le C\,.
\end{equation}
Then there exists $N>0$ such that $\{Z_n\}_{|n|\ge N}$
is a Riesz sequence in $L^2(0,T)$.

We prove that $\{Z_n\}_{n\in\mathbb{Z}'}$ is $\omega$-independent i.e. we prove that $\{\alpha_n\}_{n\in\mathbb{Z}'}=0$ when  $\{\alpha_n\}\in l^2(\mathbb{Z}')$ and

\begin{equation}\label{serie}
\sum_{n\in\mathbb{Z}'}\alpha_nZ_n=0\quad \mbox{i.e.}\quad 
\sum_{n\in\mathbb{Z}'}\alpha_n \left( z_n+ \frac{i}{n}z'_n\right)=0 
\end{equation}
 (this proof can be found in~\cite{Pandcina} althout not in such direct form, and it is reported for completeness).

Using $T>2\pi$ 
and ~\cite[Lemma~3.4]{Pandolfi} applied twice
 it is possible to prove that  $\alpha_n=\frac{\gamma_n}{n^2}$ with $\{\gamma_n\}\in l^2(\mathbb{Z}')$ 
 (see also~\cite{PandIEOT}).    
This fact justifies the termwise differentiation of the series~(\ref{serie}). Using
\begin{equation}\label{eq:DizNsecondORD}
z_n''(t)=-n^2 N(t)- n^2\int_0^t N(t-s) z_n'(s)\ZD s
\end{equation}
we get

\begin{equation} \label{serie2}
\int_0^t N(t-s)\left [\sum_{n\in\mathbb{Z}'} 
 \gamma_n\left ( z_n(s)+\frac{i}{n} z_n'(s)\right ) 
\right ]\ZD s  
-iN(t)\left [\sum_{n\in\mathbb{Z}'}  \frac{\gamma_n}{n}
\right ]
=0\,. 
\end{equation}
Computing with $t=0$ we see that
$
\sum_{n\in \mathbb{Z}'} n\alpha_n=\sum_{n\in \mathbb{Z}'}\frac{\gamma_n}{n }=0 
$ and so,
using $N(0)\neq 0$,  we get
\[
\sum_{n\in\mathbb{Z}'}\left [n^2\alpha_nz_n(s)+in\alpha_nz'_n(s)\right ]=0\quad \mbox{hence } \ 
 \sum_{n\neq \pm1,\ n\in\mathbb{Z}'}\alpha_n(n^2-1)\left [z_n +  \frac{iz_n'}{n}\right ]=0 \,.
\]
Note that $\{\alpha_n(n^2-1)\}=\{\alpha^{(1)}_n\}\in l^2(\mathbb{Z}')$.  Hence we can start a boostrap argument and repeat this procedure. After at most $2N$ iteration of the process we get
\[\sum_{|n|>N} \alpha^{(N)}_nZ_n=0\]
and so  $\alpha^{(N)}_n=0$ when $|n|>N$ since we noted that $\{Z_n\}_{|n|>N}$ is a Riesz sequence in $L^2(0,T)$. We have $\alpha_n^{(N)}=0$ if and only if $ \alpha_n=0$ and this shows that the series~(\ref{serie}) is a finite sum,
$
\sum_{n\in\mathbb{Z}',\ |n|\le N}\alpha_nZ_n=0 
$.
The proof is now finished since  it is easy to prove, as in~\cite{Pandolfi,PandolfiSHARP}, that \emph{the sequence $\{Z_n(t)\}_{n\in\mathbb{Z}'}$ is linearly independent.}

\paragraph{Step 2: $\{Z_n\}_{n\in\mathbb{Z}}$ is a Riesz sequence}
Of course,  
$\{ Z_n\}_{n\in\mathbb{Z}}$ is quadratically close to $\{ e^{\gamma t}e^{int}\}_{n\in \mathbb{Z}}$. It remains to prove $\omega$-independence, when $T>2\pi$. We prove $\{\alpha_n\} _{n\in\mathbb{Z}}=0$ when $\{\alpha_n\}\in l^2(\mathbb{Z})$ and 

\begin{equation}\label{omegind0}
\alpha_0+\sum_{n\in\mathbb{Z}'}\alpha_nZ_n=0\,.
\end{equation}
   Using that constant functions belong to $H^1$ and~\cite[Lemma~3.4]{Pandolfi} applied twice we see that $\alpha_n=\gamma_n/n^2$, $\{\gamma_n\}\in l^2$.  So, we can compute termwise the derivatives of both the sides of~(\ref{omegind0}) and we get

\begin{equation}\label{omegind1}  
\begin{split}
\sum_{n\in\mathbb{Z}'}\alpha_n\left(z'_n(t)+\frac{i}{n}\left[-n^2N(t)-n^2\int_0^tN(t-s)z_n'(s)\ZD s\right]\right)=0\, .
\end{split}
\end{equation}
Computing with $t=0$     we get
$\sum_{n\in\mathbb{Z}'}\alpha_nn=0$. 
Then (using~(\ref{eq:dizN})) the equation (\ref{omegind1}) becomes  
\begin{equation*}
\int_0^tN(t-s)\left[\sum_{n\in\mathbb{Z}'}\left (\alpha_nn^2z_n(s)+i  \alpha_n  nz'_n(s)\right )\right]\ZD s=0 
\end{equation*}
so that (using again $N(0)\neq 0$ and $\{\alpha_nn^2\}\in l^2$)

\begin{equation}\label{omegind2}
 \sum_{n\in\mathbb{Z}'}\alpha_nn^2\left [ z_n(t)+i\frac{1}{n}z'_n(t)\right ]=\sum _{n\in\mathbb{Z}'}\alpha_n n^2 Z_n(t)=0\,.
\end{equation}
The fact that $\{Z_n(t)\} _{n\in\mathbb{Z}'}$ is a Riesz sequence implies that $\{\alpha_n\}=0$ and so also $\alpha_0=0$, as we wanted to prove.

This ends the proof of Theorem~\ref{teo:regularity}.

\medskip

\noindent{\bf Remark:}
{\rm
The proof we have given shows that the cascade system
\[
\begin{cases}
w''(x,t)=w_{xx}(x,t)+\int_0^tM(t-s)w_{xx}(x,s)\ZD s\,,\quad w(0,t)=y(t)\,, \quad w(\pi,t)=0 \\
y'(t)=g(t)
w(x,0)=0\,, \quad w'(x,0)=0\,.\\
\end{cases}
\]
is controllable in $H^1_0(0,\pi)\times L^2(0,\pi)$, using the controls $g\in L_0(0,T)$ where
\[
L_0(0,T)=\left \{ g\in L^2(0,T)\,,\quad \int_0^T g(s)\ZD s=0\right \}\,.
\]
We refer to~\cite[Ch.~11]{Tucsnak} for the analysis of controllability of interconnected systems.

}

\section{Acknowledgment}

 This papers fits into the research program of the GNAMPA-INDAM and has been written in the framework of the   ``Groupement de Recherche en Contr\^ole des EDP entre la France et l'Italie (CONEDP-CNRS)''.

\end{document}